\documentclass[]{amsart}
\usepackage{geometry}                
\usepackage{graphicx}
\usepackage{epstopdf}
\usepackage{tikz}
\usepackage{pgfplots}
\usepgfplotslibrary{groupplots}
\usepackage{subcaption}
\usepackage{multirow}
\usepackage{hyperref}
\usepackage{cleveref}
\usepackage{todonotes}
\usepackage{subcaption}

 \crefname{equation}{Equation}{Equations}
 \crefname{theorem}{Theorem}{Theorems}
 \crefname{proposition}{Proposition}{Propositions}
 \crefname{problem}{Problem}{Problems}
 \crefname{lemma}{Lemma}{Lemmas}
 \crefname{definition}{Definition}{Definitions}
 \crefname{corollary}{Corollary}{Corollaries}
 \crefname{table}{Table}{Tables}
 \crefname{figure}{Figure}{Figures}

\graphicspath{ {./Figs/}}

\newcommand{\calR}{\mathcal{R}}
\newcommand{\calF}{\mathcal{F}}
\newcommand{\nuF}{\mathcal{F}_{\text{NU}}}
\newcommand{\argmin}{\operatornamewithlimits{arg\,min}}

\title{Accelerating optimization-based computed tomography via sparse matrix approximations}	\thanks{This manuscript has been authored by UT-Battelle, LLC under Contract No. DE-AC05-00OR22725 with the U.S. Department of Energy. The United States Government retains and the publisher, by accepting the article for publication, acknowledges that the United States Government retains a non-exclusive, paid-up, irrevocable, worldwide license to publish or reproduce the published form of this manuscript, or allow others to do so, for United States Government purposes. The Department of Energy will provide public access to these results of federally sponsored research in accordance with the DOE Public Access Plan (http://energy.gov/downloads/doe-public-access-plan). }
\author{Richard C Barnard \and Rick Archibald}
	\address{Computer Science and Mathematics Division, Oak Ridge National Laboratory, One Bethel Valley Road, P.O. Box 2008, MS-6211, Oak Ridge, TN 37831-6211.}
\email{barnardrc@ornl.gov} \email{archibaldrk@ornl.gov}
\begin{document}
\maketitle
\begin{abstract}
	Variational formulations of reconstruction in computed tomography have the notable drawback of requiring repeated evaluations of both the forward Radon transform and either its adjoint or an approximate inverse transform which are relatively expensive. We look at two methods for reducing the effect of this resulting computational bottleneck via approximating the transform evaluation with sparse matrix multiplications.  The first method is applicable for general iterative optimization algorithms.  The second is applicable in error-forgetting algorithms such as split Bregman.  We demonstrate these approximations significantly reduce the needed computational time needed for the iterative algorithms needed to solve the reconstruction problem while still providing good reconstructions.
\end{abstract}

\section{Introduction}\label{sec:intro}
In the present work, we focus on computed tomographic reconstruction methods based around solving a minimization problem of the general form
\begin{equation}
\label{eq:objmin}
\min_\mu\frac{\alpha}{2}\|M\calR\mu-f\|^2_2+\frac{1}{p}\|K\mu\|^p_p
\end{equation}
for $p=1,2$ and some linear operator $K$.
For simplicity, we focus on the case of parallel beam tomography.
Here, $\calR$ denotes the Radon transform and $M$ is a linear operator, similar to a ``row selector'', which restricts to those angles along which measurements are taken in constructing $f$, which is the given data.

Standard methods for tomographic reconstruction typically use filtered back projections \cite{Nat86}.
In \cite{AndCarNik16,And05}, the computational costs of such methods were reduced via efficient computations of the Radon transform and its inverse.
However, in the case of situations where the angular measurements are under resolved, filtered back projections lead to high frequency artifacts in the reconstructions.
Solving \cref{eq:objmin}, in the specific case where $p=1$ and $K=\nabla$, provides improved reconstructions without these artifacts, as demonstrated in \cite{BarTooNaf17,MelWahBab16}.
Such methods have the drawback of having an increased computational cost.
One source of this increase is in the need for repeatedly computing the Radon transform and its adjoint.
After describing in some detail this computational bottleneck in \cref{sec:bottleneck}, we describe two methods for reducing this cost in \cref{sec:overcome}.
The first relies on the fixed geometries of the tomographic reconstruction problem to replace relatively expensive transforms with sparse matrix multiplications. The second method focuses on taking advantage of the error-forgetting nature of algorithms such as the split Bregman \cite{GolOsh09} algorithm to reduce the computational costs even further.
After describing these methods, we
present in \cref{sec:examples} the ability of these two methods to reduce the computational costs of solving the total variation-penalized reconstruction problem.



\section{A computational bottleneck in optimization-based reconstruction}\label{sec:bottleneck}
We now turn to specifying the computational bottleneck we are interested in addressing.
First, if $p=2,$ the solution to \eqref{eq:objmin} solves the linear operator equation
\begin{equation}
	\label{eq:KKT}
	\big(\alpha\calR^*M^*M\calR+ K^*K\big)\mu=\calR^*M^*f.
\end{equation}
If $p=1,$ typically an iterative convex optimization algorithm is needed. Such algorithms often involve solving a series of smooth optimization problems.
For instance, the primal-dual method of Chambolle and Pock \cite{ChaPoc11} at each iterate will require us to solve smooth problems of the form
\begin{equation*}
	\argmin_u\frac{\|u-\tilde{u}\|^2}{2\tau}+\frac{\alpha}{2}\|M\calR u-g\|^2_2
\end{equation*}
which means solving
\begin{equation}
	\label{eq:PD_resolve}
	\big(I+\tau\alpha\calR^*M^*M\calR\big)u=\tilde{u}+\tau\alpha\calR^*M^*g.
\end{equation}
Similarly, the split Bregman method of Goldstein and Osher \cite{GolOsh09} (used specifically for neutron tomography in \cite{BarTooNaf17})  requires repeatedly solving a linear system of the form
 \begin{align}
	 \label{eq:SB_solve}
	 \big(\alpha\calR^*M^*M\calR&-\lambda K^*K\big)\mu =
	 \alpha\calR^*M^*f+\lambda K^*v.
 \end{align}

In contrast to the cases considered in \cite{GolOsh09} and \cite{ChaPoc11}, the solutions to \eqref{eq:KKT},\eqref{eq:PD_resolve}, and \eqref{eq:SB_solve} do not have closed form solutions.
This means one must resort to iterative linear solvers.
 Due to the expense of evaluating $\calR^*M^*M\calR$, a computational bottleneck arises in the iterative linear solvers.
 In the case of the primal-dual and the split Bregman methods on large-scale problems, this operator may need to be evaluated hundreds, or even thousands of times, necessitating an efficient implementation of its evaluation.

Let $M_{\theta_i}$ denote the restriction associated with a single direction $\theta_i$.
Via the Fourier slice theorem, for each angle $\theta_i$, the one dimensional Fourier transform of $M_{\theta_i}\calR\mu$ may be associated with a two dimensional Fourier transform of $\mu$ at points on a slice in the Fourier domain \cite{BraManBro00}.
Repeating for all angles used in constructing $M$ in \cref{eq:objmin}, we may
rewrite \cref{eq:objmin} as
\begin{equation}\label{eq:four_obj}
	\min_\mu\frac{\alpha}{2}\|\calF^*\mu-\hat{f}\|^2_2+\frac{1}{p}\|K\mu\|^p_p.
\end{equation}
However, in the discrete setting, $\calF_D^*\mu-\hat{f}$ is evaluated on the slices corresponding to the discrete angles for which measurements have been taken, where $\calF_D$ denotes a standard FFT. This results in $\calF_D^*\mu$ being needed on a non-Cartesian grid; we denote by $\nuF$ the associated discrete Fourier transform on nonuniform grids.
This means that the split Bregman method, for instance, requires solving the linear system
\begin{equation}\label{eq:SBFFT}
	(\alpha\nuF\nuF^*-\lambda K^*K\big)\mu =
	\alpha\nuF\hat{f}+\lambda K^*v.
\end{equation}
Similar statements hold for the analogues of \cref{eq:KKT,eq:PD_resolve}.
Therefore, this modified minimization problem still poses the same computational difficulties as mentioned above, however now the computational bottleneck lies in the evaluation of $\nuF\nuF^*$.

\section{Overcoming the computational bottleneck}\label{sec:overcome}
In this section we look at two methods for reducing the computational bottleneck outlined in \cref{sec:bottleneck}.
The first method involves replacing $\nuF^*\nuF$ with sparse matrix multiplications and standard fast Fourier transforms.
This method is applicable to speeding up the conjugate gradient methods needed in solving \cref{eq:four_obj} via any of the already mentioned iterative optimization methods.
The second method is specifically applicable to the split Bregman method.
Relying on the error-forgetting properties of the algorithm, we see dramatic reductions on the computational times needed without a loss in the overall performance of the algorithm.
In each case, we replace evaluating much of the effort  of evaluating $\nuF\nuF^*$ with sparse matrix multiplications.
From an application perspective, these matrices depend primarily on the resolution of the instruments, the frequency of angular measurements, and parameters which typically take only a few discrete values.
This means for a particular tomographic reconstruction workflow, one may precompute the relevant matrices offline, enabling significant speedup in situations where multiple reconstructions are needed.

\subsection{Implementation of convolution fusing}\label{sec:convfuse}

In \cite{GreLee04}, nonuniform fast Fourier transforms (NUFFT) were developed, providing a computationally efficient method of computing these operators. We briefly review this construction.  The one-dimensional, type 1 NUFFT of $f(x)$, is computed by taking
\begin{equation}
	F(k) = \sqrt{\frac{\pi}{\tau}}e^{k^2\tau}\calF_D(f*g_\tau)_k
\end{equation}
where $\tau$ is a parameter determining the accuracy desired, and $g_\tau$ is a periodic heat kernel of the form
\begin{equation*}
	g_\tau(x) = \sum_{\ell=-\infty}^{+\infty}e^{-(x-2\ell\pi)^2/4\tau}.
\end{equation*}
For notational convenience, we let $A$ be the linear operator defined by $A(f)(k):=\sqrt{\frac{\pi}{\tau}}e^{k^2\tau}f_k$ and let $Bf:=g_\tau*f(x).$

We may then rewrite $\nuF\nuF^*$ as
\begin{equation}\label{eq:fusedop}
	A\mathcal{F}_DBB^*\mathcal{F}_D^*A
\end{equation}
(noting that $A$ is self-adjoint).
As noted in \cite{GreLee04}, the main computational task in computing the NUFFT is evaluating
\begin{equation*}
	Bf(\pi m/M) = \sum_{j}f_j\sum_{\ell=-\infty}^\infty e^{(-x_j-\pi m/M-2\pi\ell)^2/4\tau}.
\end{equation*}
In order to do this efficiently, a parameter $M_{sp}$ is chosen in order to achieve a desired level of accuracy in the NUFFT.
For $m\in[0,2N-1]$, we approximate $Bf$ by truncated summations, taking the nearest $2M_{sp}$ points around $x_j$.
Outside of this window, the Gaussians decay sufficiently sharply enough to be ignored.

Defining $m_j=\sup\{\tilde{m}:\pi\tilde{m}/M\leq x_j \}$ and $I_m:=\{j:|m_j-m|\leq M_{sp}\}$,
$Bf(\pi m/N)$ is approximately
\begin{equation*}
	\sum_{j\in I_m}f(x_j)e^{-(x_j-\pi m/M-2\pi\ell)^2/4\tau}.
\end{equation*}
It is noted in \cite{GreLee04}, drawing from \cite{DutRok93}, that if $M_{sp}=6,~\tau=6/N^2$, the NUFFT is accurate to 6 digits of accuracy and if $M_{sp}=12,~\tau=12/N^2$, the NUFFT is accurate to 12 digits of accuracy.

We now turn to $B^*B$, which takes the form
\begin{align*}
	BB^*f(x_j)=\sum_{m\in J_j}\sum_{k\in I_m} \Big(E_1(j,m)E_1(k,m)&(E_2(j,m))^{m_j-m}(E_2(k,m))^{m_k-m}\\
	&e^{-(\pi(m_j-m)/2N)^2/\tau}e^{-(\pi(m_k-m)/2N)^2/\tau}\Big)f_k,
\end{align*}
where $J_j = \{m:|m_j-m|\leq M_{sp}\}.$
This operator is sparse; only grid points within $M_{sp}$ of $x_j$ are used, along with appropriate weights taken from the Gaussian sources.  In higher dimensions, this sparsity structure holds; for instance in two dimensions, the only grid points used in computing $B^*Bf(x_{i,j})$ are those points such that $|m-i|<M_{sp}$ and $|j-n|<M_{sp}.$
The structure of $B^*B$ is determined by the radial resolution of the data $f$, the set of angular measurements used, and the value $M_{sp}$.
The first depends on the accuracy of the instrumentation, which typically is fixed.
In situations where the set of angular measurements is either fixed, or from one of several canonical sampling regimes, the $B^*B$ may be precomputed as a sparse matrix.
The resulting sparse matrix fuses the convolutions into a single operation.
Thus, the NUFFT's needed in solving \cref{eq:SBFFT} may be replaced by 2 FFT's,  a single sparse matrix-vector multiplication, and 2 diagonal matrix-vector multiplications.

\subsection{Acceleration via surrogate matrices}
Taking advantage of fusing the convolutions as outline above in \cref{sec:convfuse} means, at each iterate in the split Bregman algorithm, we now solve
\begin{equation}\label{eq:SB_update}
	\big(\alpha A\mathcal{F}_D(B^*B)\mathcal{F}_D^*A-\lambda K^*K\big)\mu = \alpha \nuF^* \hat{f}+\lambda K^*v,
\end{equation}
where now $B^*B$ has been replaced with a sparse matrix multiplication.
However, the (split) Bregman algorithm has an error forgetting property \cite{GolOsh09}; that is, one only needs in general to solve \cref{eq:SB_update} approximately in order to still have a convergent algorithm.
The benefit for doing exact Bregman updates lies in the ability of Bregman to converge in two exact updates, for sufficiently small initial Bregman distance, as discussed in \cite{YinOsh13}.
However, we note that solving \cref{eq:SB_update} exactly (to machine precision) can be quite expensive.
For large-scale imaging problems, the expense of exactly solving \cref{eq:SB_update} may offset the near immediate convergence of Bregman updates for good (with respect to the Bregman distance) initialization.
Indeed, in \cite{BarTooNaf17}, inexact Bregman updates were sufficient for a similar formulation of the tomographic reconstruction problem.

In order to further reduce the computational costs associated with solving \cref{eq:SB_update}, we look to leverage the error-forgetting nature of the split Bregman algorithm by precomputing a sparse matrix $T$ which approximates $A\calF_D(B^*B)\calF_D^*A$.
We first consider the cases where $u$ is the simple $N\times N$ image given by \begin{equation*}
	 u_{i,j} = 	\begin{cases}
 					1 & i=N/2,j=N/2\\
				 	0 & else
			   	\end{cases}.
\end{equation*}
The result of $\calF^*A$ being applied to this $u$ is constant function.
As described above in \cref{sec:convfuse}, $B^*B\calF^*_DA(\delta_{i,N/2}\delta_{j,N/2})$ will then be the sum of Gaussians.
\begin{figure}
	\centering
	\begin{subfigure}[t]{.30\linewidth}
		\includegraphics[width = \textwidth]{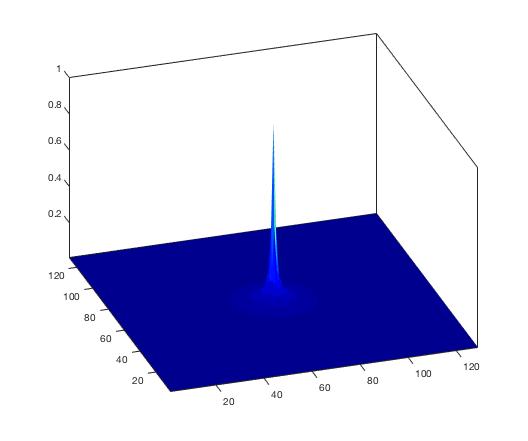}
		\caption{Closeup of center of a normalization of $B^*B\calF^*_DA u$ for $N=256$ and $M$ corresponding to 403 evenly spaced angular samples.}
		\label{fig:op_pt_gauss}
	\end{subfigure}
	\hfill
	\begin{subfigure}[t]{.30\linewidth}
		\includegraphics[width = \textwidth]{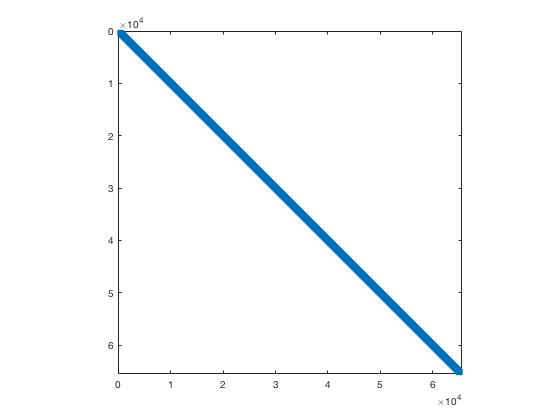}
		\caption{Sparsity pattern for $T$ ignoring periodic boundary extension.}
		\label{fig:surrogate_sparse}
	\end{subfigure}
	\hfill
	\begin{subfigure}[t]{.30\linewidth}
		\includegraphics[width = \textwidth]{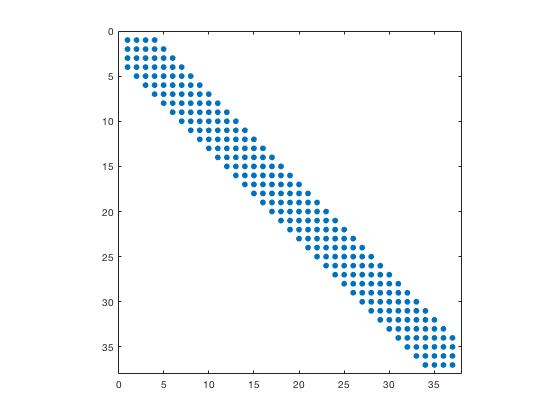}
		\caption{Closeup of sparsity pattern for $T$.}
		\label{fig:surrogate_sparse_close}
	\end{subfigure}
	\caption{The sparsity of $T$ arising from the peaked nature of $B^*B\calF^*_DA u$. }
\end{figure}
Therefore, $\big(A\calF_D(B^*B)\calF_D^*A\big)(u)$ will be a (possibly rescaled) Gaussian centered at $(N/2,N/2)$.
This is shown in \cref{fig:op_pt_gauss} for $N=256$ where we show the central $128\times128$ region when $M$ corresponds to a Shannon sampling frequency in angle.
We will use the peakedness seen in that image to construct a $T$.
First, for a chosen truncation radius $r>0,$ define for some $k,\ell\in[-N/2,N/2]$ the matrix $T^{k,\ell}$ by
\begin{equation*}
	T^{k,\ell}_{m,n}=\begin{cases}
	\big[\big(A\calF_D(B^*B)\calF_D^*A\big) (\delta_{i,N/2}\delta_{j,N/2})\big]_{m-k,n-\ell} & |(m-k,n-\ell)-(N/2,N/2)|\leq r\\
	0& \text {otherwise}
	\end{cases}.
\end{equation*}
We then define
\begin{equation}
\label{eq:def_sparse_surr}
T:=\sum_{k,\ell=-N/2}^{N/2} T^{k,\ell}
\end{equation}
as the sum of these translated stencil matrices. This $T$ is used to compute split Bregman updates; instead of \cref{eq:SB_update}, we solve (possibly approximately)
\begin{equation}\label{eq:surr_SBupd}
	(\alpha T-\lambda K^*K)\mu = \alpha \nuF^* \hat{f}+\lambda K^* v.
\end{equation}

It is important to recall that for an an image which is a single non-zero pixel near the image's edge, e.g. $\tilde{u}_{i,j}=\delta_{i,k}\delta_{j,\ell}$ for some $k<r$,
$\big(A\calF_D(B^*B)\calF_D^*A\big)\tilde{u}$ will have nontrivial entries near the opposite edge, due to the periodicity of the operators.
Our choice of $T$ suppresses these values on the opposite edge.
The motivation for this choice is to induce a banded structure in $T$, as seen in \cref{fig:surrogate_sparse,fig:surrogate_sparse_close}; in this case $r=3,$ resulting in a banded matrix with 7 diagonal and off-diagonal entries per row.
However, this is reasonable in the case of computed tomography where the sample to be measured is fully contained in the interior of the image.
The choice of suppressing values on the opposite edge corresponds to an assumption that the object of interest in the reconstructed image is at least $r$ pixels from each edge of the image and that near the edge, the true $\mu$ can be assumed to be $0$.

\section{Numerical Implementation}\label{sec:examples}
In order to test the efficiency of using the sparse linear matrix representation, we solve \cref{eq:objmin} for $K,p$ corresponding to anisotropic TV regularization via a split-Bregman algorithm (as discussed for the general problem in \cite{GolOsh09} and specifically for computed tomography in \cite{VanGooBee11,BarTooNaf17}). Throughout, we use the Shepp-Logan phantom as the reference data and $\alpha=\lambda=1$.  All computations reported below were performed in MATLAB on a workstation using Intel 2.20 GHz Xeon Processors with 55 MB caches.

\begin{figure}
	\centering
	\input{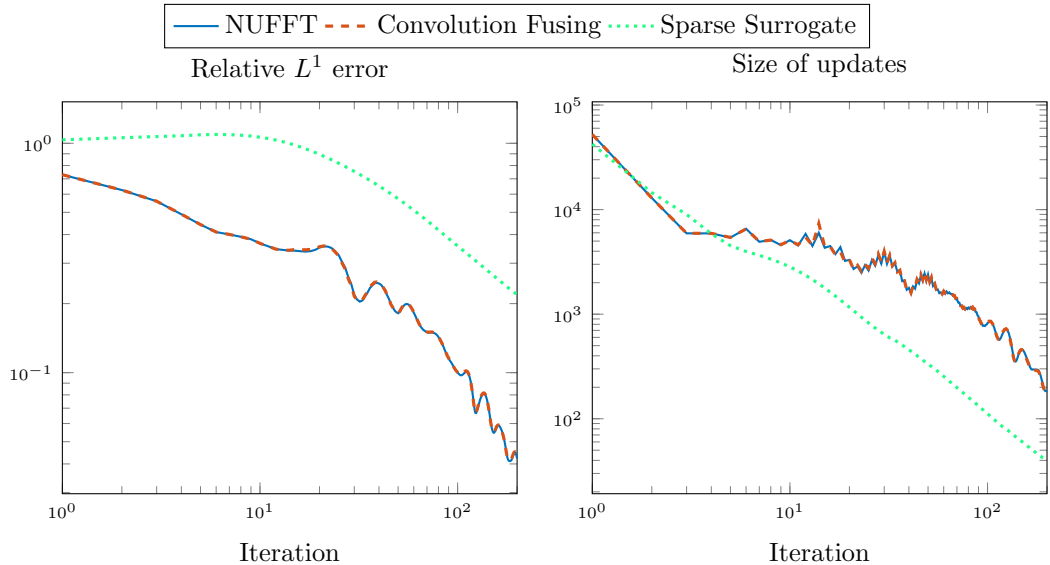}
	\caption{Behavior of algorithm over iterations during Split Bregman for each method of solving \eqref{eq:SBFFT}, $N_x = 512$ and $6$ digits of accuracy for NUFFT.  On the left is the $L^1$ error at each update; on the right, the $L^1$ norm of the Bregman updates.}
	\label{fig:fixed_update_comp}
\end{figure}

First, in order to compare the cost and accuracy of using the NUFFT as opposed to either the convolution fusing method or the use of sparse surrogates,  we look at the time required to compute a fixed number of inexact Bregman iterations.  We compute 200 inexact Bregman iterations; each of these iterations involves approximating the solution to \cref{eq:SB_solve} via 5 conjugate gradient steps.
Therefore, 1000 evaluations of $\alpha\calR^*M^*M\calR$ are needed.
We perform this for various image sizes as well as various values of $M_{sp}$ used in computing the NUFFT and its corresponding $B$.
For each value of $N_x,$ we choose $M$ that corresponds to using $N_x\frac{\pi}{4}$ evenly spaced angular measurements.
We summarize the results in \cref{tab:fixed_iter_time,tab:fixed_iter_error} as well as \cref{fig:fixed_update_comp}.
First, in \cref{tab:fixed_iter_time}, we report the time required for performing 200 inexact Bregman iterations as well the portion of that time required for performing the conjugate gradient computations.
As expected, increasing $M_{sp}$ increases the time required when using the NUFFT and convolution fusing-based optimizations.  The sparsity structure of $T$ in the surrogate-based method is dominated by the choice of $r$ (here set to be $1$), resulting in the time needed being unchanged for different $M_{sp}$.
Additionally, the time required per evaluation of $\alpha\calR^*M^*M\calR$ using the NUFFT is, for $M_{sp}=2$, roughly 3 or 4 times longer than using convolution fusing and 10-20 times longer than using sparse surrogates.
At higher values of $M_{sp},$ the speedup using convolution fusing seems to hold; meanwhile the possible speedup by using sparse surrogates becomes one of several orders.
{
\begin{table}[htbp]
	\caption{Total time needed to complete 200 Split Bregman iterations and total time needed during conjugate gradient steps only }
	\label{tab:fixed_iter_time}%
	{{\footnotesize\renewcommand{\arraystretch}{1.1}
		\centering
\begin{tabular}{lrrrrrr}
	\hline
	Method  & \multicolumn{2}{c}{$N_x = 128$} & \multicolumn{2}{c}{$N_x = 256$} & \multicolumn{2}{c}{$N_x = 512$} \\
	\hline
	& Total  & CG  & Total & CG  & Total  & CG \\
	\hline
	NUFFT, $M_{sp}=2$ &  33.2305 & 27.9056 & 142.0550 & 118.9043 & 581.6318 & 479.4620 \\
	Convolution Fusing, $M_{sp}=2$ & 9.6462 & 7.7387 & 43.5623 & 34.7016 & 199.0465 & 154.3572 \\
	Sparse Surrogate, $M_{sp}=2$ & 3.6270 & 1.8375 & 18.1970 & 9.1952 & 93.2822 & 49.6176 \\
	\hline
	NUFFT, $M_{sp}=6$ & 115.6706 & 98.5060 & 462.2991 & 393.7125 & 1810.3551 & 1530.83 \\
	Convolution Fusing, $M_{sp}=6$ & 37.3927 & 34.1586 & 160.1846 & 145.3728 & 666.03458 & 598.0480 \\
	Sparse Surrogate, $M_{sp}=6$ & 4.8572 & 1.7213 & 24.7311 & 9.0292 & 116.5725 & 48.3647 \\
	\hline
	NUFFT, $M_{sp}=12$ & 371.3576 & 317.8252 & 1470.7238 & 1258.3341 & 6334.8399 & 5403.2094 \\
	Convolution Fusing, $M_{sp}=12$& 130.4410 & 122.7075 & 658.4656 & 619.9470 & 3122.7715& 2951.3699 \\
	Sparse Surrogate, $M_{sp}=12$ & 9.8901 & 1.7619 & 39.4063 & 7.9888 & 204.9899 & 47.2547 \\
	\hline
\end{tabular}
}
}%

\end{table}%

}
{
\begin{table}[htbp]
	\centering
\caption{Relative $L^1$ errors after 200 Split Bregman iterations, 5 CG steps per iteration}
\label{tab:fixed_iter_error}%
{\renewcommand{\arraystretch}{1.1}
\begin{tabular}{lrrr}
	\hline
	Method  & \multicolumn{1}{l}{$N_x = 128$} & \multicolumn{1}{l}{$N_x = 256$} & \multicolumn{1}{l}{$N_x = 512$} \\
	\hline
	NUFFT, $M_{sp}=2$ & 0.090699645 & 0.056697226 & 0.042466068 \\
	Convolution Fusing, $M_{sp}=2$ & 0.090704816 & 0.056725578 & 0.042493176 \\
	Sparse Surrogate, $M_{sp}=2$ & 0.305754799 & 0.199389404 & 0.133857662 \\
	\hline
	NUFFT, $M_{sp}=6$ & 0.090241338 & 0.056122227 & 0.042219283 \\
	Convolution Fusing, $M_{sp}=6$ & 0.090235638 & 0.056130045 & 0.042173041 \\
	Sparse Surrogate, $M_{sp}=6$ & 0.304973205 & 0.199030832 & 0.133646147 \\
	\hline
	NUFFT, $M_{sp}=12$ & 0.090241343 & 0.056122246 & 0.042219412 \\
	Convolution Fusing, $M_{sp}=12$ & 0.090235639 & 0.056130049 & 0.042173149 \\
	Sparse Surrogate, $M_{sp}=12$ & 0.304972812 & 0.199030794 & 0.133646149 \\
	\hline
\end{tabular}%
}
\end{table}%

}

However, the speedup using sparse surrogates is tempered by a significant loss of accuracy in the reconstruction.
This can be seen in \cref{tab:fixed_iter_error}, where we report the relative $L^1$ error of the reconstruction using these 200 split Bregman iterations.
The sparse surrogate-based reconstruction results in a relative error that is between 2 and 5 times higher than reconstructions using the other methods.
One perhaps surprising result from this test is that the choice of $M_{sp}$ appears to not significantly affect the quality of the reconstruction after a fixed number of iterations.
Finally, in \cref{fig:fixed_update_comp}, we display the relative errors and the $L^1$ norm of the updates over the iterations.
Using either the NUFFT or convolution fusing results in nearly the same iterative reconstruction.
Additionally, despite the higher error in the sparse surrogate method, the convergence rates are unchanged.  Indeed, the sparse surrogate method leads, after a handful of iterations, to a iteration scheme which is largely monotonically reducing the error of the reconstruction, unlike the other two methods.  These observations suggest that for a fixed number of iterations, convolution fusing provides a significant speedup over the NUFFT while not impacting the performance of the image reconstruction procedure.
Meanwhile, the use of sparse surrogates decreases the accuracy at any given iterate.
This also holds if one increases $r$: for $r=3,$ $200$ Bregman updates require $136.0654$ seconds results in a relative $L^1$ error of $0.129994$ error.
As expected, then, increasing $r$ results in more expensive (due to decreasing the sparsity of the resulting $T$) albeit more accurate iterations. 
However, due to the shorter time of each iteration in the sparse surrogate method, there is a possibility that the overall performance of solving \eqref{eq:objmin} can be improved: we may replace relatively accurate iterations with more iterations which, individually, are significantly quicker to perform.

To test this, we look to solve \eqref{eq:four_obj} for $512\times 512$ images and $K,p$ again corresponding to anisotropic total variation regularization.
We consider $M$ corresponding to taking $512\frac{\pi}{4d}$ angular measurements for $d=1,2,4$ and $T$ corresponding to $r=3.$
That is, we consider reconstruction problems where the angular measurements are under-resolved to differing degrees.
In \cite{BarTooNaf17}, the highly under-resolved case was of interest, as angular measurements are relatively expensive in some applications.
We will compare the results of using split Bregman to solve \eqref{eq:four_obj} with both convolution fusing and with sparse surrogates.
That is, we will compare using \eqref{eq:SB_update} with using \eqref{eq:surr_SBupd} to compute the Bregman update; again, we use at most 5 conjugate gradient steps to compute this update.
As a stopping criterion, we terminate either after a maximum of 6000 Bregman updates or if the $L^1$ norm of a Bregman update is less than $10^{-8}$ of the first Bregman update.
In \cref{tab:Sub1_cut_tab,tab:Sub2_cut_tab,tab:Sub4_cut_tab} we list the results of using each method to compute the Bregman update for $d=1,2,4$, respectively.
In each we show the time and number of iterations needed for Bregman updates to fall below a certain threshold, as well as the relative $L^1$ error of the reconstruction at that iteration. Additionally we show in \cref{fig:full_update_comp,fig:full_sub2_update_comp,fig:full_sub4_update_comp} the $L^1$ error and size of the Bregman updates as a function of computational time for $d=1,2,4$, respectively.

{
\begin{table}[tbp]
	\caption{Time and iterations needed for Bregman updates to fall below listed tolerances, when $N_x=512$ and $N_x=512$ and $N_x\frac{\pi}{4}$ angular measurements.   $*$ denotes failure of convergence to listed tolerance, along with time needed to reach maximum iteration. }
	\label{tab:Sub1_cut_tab}%
	{\small \renewcommand{\arraystretch}{1.1}
  \centering
    \begin{tabular}{llrrrr}
    	\hline
    	Method & & \multicolumn{4}{c}{Bregman update size}\\
    	\hline
	 &   & \multicolumn{1}{c}{$10^{-2}$} & \multicolumn{1}{c}{$10^{-4}$} & \multicolumn{1}{c}{$10^{-6}$} & \multicolumn{1}{c}{$10^{-8}$} \\
    \multirow{3}{*}{Convolution Fusing} & Iterations needed & 676   & 5984 & (*)6000 & (*) \\
    & Time Needed & 2368.153445 & 22560.0831 & 22714.6989 & (*) \\
    & Rel. $L^1$ error & 0.035400827 & 
    0.00064873
     & 0.00064175 & (*) \\
    \hline
    \multirow{3}{*}{Sparse Surrogate} & Iterations needed & 65    & 726   & 1835  & 3365 \\
    & Time Needed & 44.807082 & 489.120469 & 1199.41347 & 2175.546164 \\
    & Rel. $L^1$ error & 0.271823575 & 0.0189007 & 8.47635E-05 & $1.6814\times 10^{-5}$ \\
    \hline
    \end{tabular}%
}

\end{table}%

}
{
\begin{table}[tbp]
	\caption{Time and iterations needed for Bregman updates to fall below listed tolerances, when $N_x=512$, with $N_x\frac{\pi}{8}$ angular measurements.   $(*)$ denotes failure of convergence to listed tolerance, along with time needed to reach maximum iteration. }
	\label{tab:Sub2_cut_tab}
	{\footnotesize \renewcommand{\arraystretch}{1.1}
	\centering
	\begin{tabular}{llrrrr}
		\hline
		Method & & \multicolumn{4}{c}{Bregman update size}\\
		\hline
		&   & \multicolumn{1}{c}{$10^{-2}$} & \multicolumn{1}{c}{$10^{-4}$} & \multicolumn{1}{c}{$10^{-6}$} & \multicolumn{1}{c}{$10^{-8}$} \\
		\multirow{3}{*}{Convolution Fusing} & Iterations needed & 60    & 589   & $(*)$6000 & $(*)$\\
	    & Time Needed & 221.73419 & 2071.667173 & 18517.54548 & $(*)$ \\
	    & Rel. $L^1$ error & 0.259376623 & 0.00648995 & $4.23422\times10^{-5}$ & $(*)$ \\
	    \hline
	    \multirow{3}{*}{Sparse Surrogate} & Iterations needed & 19    & 270   & 870   & 1765 \\
	    & Time Needed & 12.094469 & 168.529048 & 538.209406 & 1114.983433 \\
	    & Rel. $L^1$ error & 0.409244885 & 0.034075116 & 0.000182022 & $1.7571\times 10^{-6}$ \\
	   \hline
    \end{tabular}%
}
\end{table}%

}
{
\begin{table}[tbp]
	\caption{Time and iterations needed for Bregman updates to fall below listed tolerances, when $N_x=512$, with $N_x\frac{\pi}{16}$ angular measurements.   $(*)$ denotes failure of convergence to listed tolerance, along with time needed to reach maximum iteration. }
	\label{tab:Sub4_cut_tab}
	{\small \renewcommand{\arraystretch}{1.1}
	\centering
	\begin{tabular}{llrrrr}
		\hline
		Method & & \multicolumn{4}{c}{Bregman update size}\\
		\hline
		&   & \multicolumn{1}{c}{$10^{-2}$} & \multicolumn{1}{c}{$10^{-4}$} & \multicolumn{1}{c}{$10^{-6}$} & \multicolumn{1}{c}{$10^{-8}$} \\
		\multirow{3}{*}{Convolution Fusing} & Iterations needed & 7     & 159   & 1116  & 4815 \\
	    & Time Needed & 20.335822 & 456.510441 & 2926.561617 & 11480.16517 \\
	    & Rel. $L^1$ error & 0.589578325 & 0.019544919 & 0.000134007 & $2.96197\times 10^{-6}$ \\
	    \hline
	    \multirow{3}{*}{Sparse Surrogate} &  Iterations needed & 3     & 86    & 471   & 1195 \\
	    &  Time Needed & 1.821839 & 49.449029 & 270.114974 & 682.877529 \\
	    & Rel. $L^1$ error & 0.792041713 & 0.070293249 & 0.001610545 & $3.48151 \times 10^{-5}$ \\
	    \hline
    \end{tabular}%
}

\end{table}%

}

\begin{figure}
	\centering
	\input{./Figs/Combo_time_Upd_Sub1.tex}
	\caption{Error and size of Bregman updates as functions of needed computational time, with $N_x=512$ and $N_x\frac{\pi}{4}$ angular measurements.  On the left is the relative $L^1$ error at each update; on the right, the $L^1$ norm of the Bregman updates.}
	\label{fig:full_update_comp}
\end{figure}
\begin{figure}
	\centering
	\input{./Figs/Combo_time_Upd_Sub2.tex}
	\caption{Error and size of Bregman updates as functions of needed computational time, with $N_x=512$ and $N_x\frac{\pi}{8}$ angular measurements.  On the left is the relative $L^1$ error at each update; on the right, the $L^1$ norm of the Bregman updates.}
	\label{fig:full_sub2_update_comp}
\end{figure}
\begin{figure}
	\centering
	\input{./Figs/Combo_time_Upd_Sub4.tex}
	\caption{Error and size of Bregman updates as functions of needed computational time, with $N_x=512$ and $N_x\frac{\pi}{16}$ angular measurements.  On the left is the relative $L^1$ error at each update; on the right, the $L^1$ norm of the Bregman updates.}
	\label{fig:full_sub4_update_comp}
\end{figure}

We see that for any fixed time budget, the sparse surrogate-based method provides improved accuracy in very few iterations.
In the case of $d=1,$ the time needed to get $10^{-5}$ order relative $L^1$ error via the sparse surrogate method is approximately $10\%$ of the time needed to reach $10^{-4}$ $L^1$ error via convolution fusing.
Similar observations hold for $d=2.$  However, for $d=4,$ the sparse surrogate-based method, while taking vastly less time than the convolution fusing-based method, results in a final result which is more inaccurate.
This is due to prematurely terminating the reconstruction method, as the optimization routine detects a convergence due to the relatively small size of the Bregman updates.
Nonetheless, the time needed to reach $3.4739\times 10^{-5}$ error via convolution fusing is  $5104$ seconds, versus the $683$ seconds needed for sparse surrogates.

\section{Conclusions}
Approximating the transform $\calF^*\calF$ with sparse matrices reduces the computational bottleneck that arises in solving variational tomographic reconstruction problems.
Convolution fusing replaces the most computationally expensive portion of the evaluation of this operator with a sparse matrix multiplication, reducing the computation time without any additional loss of accuracy.
In cases where an error-forgetting algorithm is to be used, sparse surrogates further reduce the cost of the computation by replacing the entire operator with a sparse matrix; this comes at the cost of a loss of accuracy.
However, the reduced time needed per iteration offsets the need for additional iterations, especially in the case of only moderately unresolved angular measurements.
In the case of low angular resolution reconstructions where the spares surrogates result in prematurely terminating algorithms, it may prove useful as a warm-start strategy.
That is, as we expect sparse surrogate-based iterations may prematurely converge, we use the result as an initial guess for fused convolution-based iterations which are more expensive.


\bibliographystyle{plain}
\bibliography{fused}

\end{document}